\documentclass{aims}
\usepackage{amsmath}
 \usepackage{amscd}
 \usepackage{bbm}
  \usepackage{paralist}
  \usepackage{graphics} 
  \usepackage{epsfig} 
\usepackage{graphicx}  \usepackage{epstopdf}
 \usepackage[colorlinks=true]{hyperref}
\hypersetup{urlcolor=blue, citecolor=red}

\def\currentvolume{X}
 \def\currentissue{X}
  \def\currentyear{200X}
   \def\currentmonth{XX}
    \def\ppages{X--XX}
     \def\DOI{10.3934/xx.xx.xx.xx}

\newtheorem{theorem}{Theorem}[section]

\newtheorem{q}{Question}
\theoremstyle{definition}

\newcommand{\ep}{\varepsilon}

\newcommand{\R}{\mathbb R}

\title[nonlocal-to-local limit] 
      {Recent results on the singular local limit for nonlocal conservation laws}

\author[Maria Colombo, Gianluca Crippa, Marie Graff and Laura V. Spinolo]{}

\subjclass{Primary: 35L65, 65M08.}
 \keywords{nonlocal conservation law, traffic model, downstream traffic density, numerical viscosity, 
singular limit, local limit}

 \email{maria.colombo@epfl.ch}
 \email{gianluca.crippa@unibas.ch}
 \email{marie.graff@auckland.ac.nz}
 \email{spinolo@imati.cnr.it}

\thanks{$^*$ Corresponding author}

\begin{document}
\maketitle

\centerline{\scshape Maria Colombo}
\medskip
{\footnotesize
 \centerline{EPFL SB, Station 8, CH-1015 Lausanne, Switzerland}
  }
  
\medskip

\centerline{\scshape Gianluca Crippa}
\medskip
{\footnotesize
 \centerline{Departement Mathematik und Informatik, Universit\"at Basel,}
\centerline{Spiegelgasse 1, CH-4051 Basel, Switzerland}
  }

\medskip

\centerline{\scshape Marie Graff}
\medskip
{\footnotesize
 \centerline{Department of Mathematics, University of Auckland,} 
 \centerline{Private Bag 92019, Auckland 1142, New Zealand}
  }

\medskip

\centerline{\scshape Laura V. Spinolo$^*$ }
\medskip
{\footnotesize
 \centerline{IMATI-CNR, via Ferrata 5, I-27100 Pavia, Italy}
  }

\bigskip


\begin{abstract}
We provide an informal overview of recent developments concerning the singular local limit of nonlocal conservation laws. In particular, we discuss
some counterexamples to convergence and we highlight the role of numerical viscosity in the numerical investigation of the nonlocal-to-local limit. We also state some open questions and describe recent related progress. 
\end{abstract}

\section{Introduction}
We consider the nonlocal conservation law 
\begin{equation}
\label{e:nl}
     \partial_t u + \partial_x \big[ u V({ u\ast \eta}) \big] =0.
\end{equation}
In the previous expression, the unknown is the function $u: \R_+ \times \R \to \R$, the function $V: \R \to \R$ is Lipschitz continuous and the term 
$u \ast \eta$ is the convolution, computed with respect to the space variable $x$ only, of the solution $u$ with the convolution kernel $\eta: \R \to \R$.  
For the time being we assume that $\eta$ satisfies the following assumptions:
\begin{equation}
\label{e:eta} 
   \eta \in C^1_c (\R), \quad \quad \eta \ge 0, \quad \quad \int_\R \eta(x) dx =1, 
\end{equation}
but actually the regularity assumptions on $\eta$ can be relaxed, as we will see in~\S\ref{s:ak}. 
Nonlocal equations in the form~\eqref{e:nl} have been extensively studied in recent years owing to the applications to (among others) models of sedimentation and pedestrian and vehicular traffic, see for instance~\cite{Sedimentation,BlandinGoatin,ChiarelloGoatin,ColomboGaravelloMercier,ColomboHertyMercier} and the references therein.  

Consider the Cauchy problem posed by coupling~\eqref{e:nl} with the initial datum 
\begin{equation}
\label{e:id}
      u(0, x) = \bar u (x).  
\end{equation}
Existence and uniqueness results have been obtained in various frameworks by several authors, see among others~\cite{BlandinGoatin,ColomboHertyMercier,CLM,KeimerPflug}. 

In this note we review some recent progress in the analysis of the singular local limit of~\eqref{e:nl}, which is defined as follows. Fix a parameter $\ep>0$,  consider the rescaled function
$$
     \eta_\ep (x)= \frac{1}{\ep} \eta \left(  \frac{x}{\ep}  \right) 
$$
and note that, owing to the third condition in~\eqref{e:eta}, when $\ep \to 0^+$ the family $\eta_\ep$ converges weakly$^\ast$ in the sense of measures  to the Dirac delta. By plugging $\eta_\ep$ into~\eqref{e:nl},\eqref{e:id} we arrive at the family of Cauchy problems 
\begin{equation}
\label{e:nlee}
     \left\{ 
     \begin{array}{ll}
     \partial_t u_\ep + \partial_x \big[ u_\ep V({ u_\ep \ast \eta_\ep}) \big] =0 \\
     u_\ep (0, x) = \bar u (x). \\
    \end{array}
    \right.
\end{equation}
We now consider the limit $\ep \to 0^+$: since $\eta_\ep$ converges to the Dirac delta, from the equation at the first line of~\eqref{e:nlee} we formally recover the nonlinear conservation law 
\begin{equation}
\label{e:l}
     \partial_t u + \partial_x \big[ u V(u) \big] =0. 
\end{equation}
The above derivation is completely formal, and whether or not it can be rigorously justified is the object of the following question, which was originally posed in~\cite{ACT}. 
\begin{q}\label{q:q} 
Does $u_\ep$, solution of~\eqref{e:nlee}, converge (in some suitable topology) to the entropy admissible solution of~\eqref{e:id},\eqref{e:l} as $\ep \to 0^+$? 
\end{q}
We refer to~\cite{Dafermos:book} for the definition of entropy admissible solution of~\eqref{e:id},\eqref{e:l}.
In this work we overview some recent developments concerning Question~\ref{q:q}. The exposition is organized as follows: in~\S\ref{s:nltl} we show that, notwithstanding numerical evidence suggesting the opposite, the answer to Question~\ref{q:q} is in general negative. In~\S\ref{s:ne} we discuss a possible explanation of the reason why the numerical evidence provides the wrong intuition. Finally, in~\S\ref{s:ak} we introduce Question~\ref{q:q3}, 
which is a refinement of Question~\ref{q:q} in a more specific setting motivated by the applications to vehicular traffic models. Question~\ref{q:q3} is still open, but 
recent progress has been recently achieved and we discuss it in~\S\ref{s:ak}.
\section{The nonlocal-to-local limit} \label{s:nltl}
Question~\ref{q:q} was originally motivated by numerical evidence. More precisely, in~\cite{ACT} the authors exhibit numerical experiments where the solution of the nonlocal Cauchy problem~\eqref{e:nlee} gets closer and closer to the entropy admissible solution of~\eqref{e:id},\eqref{e:l} as $\ep \to 0^+$, thus suggesting a positive answer to Question~\ref{q:q}. This was later confirmed by other numerical experiments, see for instance~\cite{BlandinGoatin}. 

Another positive partial answer to Question~\ref{q:q} is provided by~\cite[Proposition 4.1]{Zumbrun}, which loosely speaking states that the answer to Question~\ref{q:q} is positive provided that the convolution kernel $\eta$ is even (i.e. $\eta(x)= \eta(-x)$, for every $x$) and the limit solution $u$ is smooth. The rationale underpinning~\cite[Proposition 4.1]{Zumbrun} is basically the following. Assume that the initial datum $\bar u$ is smooth and say compactly supported, then there is a time interval $[0, T]$ where the entropy admissible solution of~\eqref{e:id},\eqref{e:l} is smooth, i.e. it is a classical solution. Proposition 4.1 in~\cite{Zumbrun} states that on the interval $[0, T]$ the family $u_\ep$ converges to $u$, in the uniform $C^0$ norm. 

Despite the above mentioned results, the answer to Question~\ref{q:q} is, in general, negative.  
More precisely, in~\cite{CCS} we exhibit three counterexamples that rule out the possibility that the family $u_\ep$ solving~\eqref{e:nlee} converge to 
the entropy admissible solution of~\eqref{e:id},\eqref{e:l}. The counterexamples are completely explicit and rule out not only strong convergence, but also i) weak convergence and ii) the possibility of extracting from $u_\ep$ a (strongly or weakly) converging subsequence. In one case we even manage to rule out the possibility that $u_\ep$ converges to a distributional solution  of~\eqref{e:id},\eqref{e:l}, i.e. we do not need to require that the limit $u$ is entropy admissible to rule out convergence. The counterexamples are constructed in~\cite[\S5.1,\S5.2,\S5.3]{CCS} and at the beginning of each of \S5.1, \S5.2 and \S5.3 the basic ideas underpinning the construction of the counterexample are overviewed. Loosely speaking the very basic mechanism is that  in each of the counterexamples we manage to single out a property that  i) is satisfied by the solution $u_\ep$ of~\eqref{e:nlee}, for every $\ep>0$; ii) is stable under weak or strong convergence, i.e. it passes to the weak or strong limit; iii) is \emph{not} satisfied by the entropy admissible solution of~\eqref{e:id},\eqref{e:l}.  The exact property verifying conditions i), ii) and iii)  is different in each counterexample: in the first one it is the fact that the integral over $\R_-$ is constant in time, in the second one the fact that $u_\ep$ is identically $0$ at positive values of $x$. Finally, in the third counterexample we single out a functional that is constant in time when evaluated at $u_\ep(t, \cdot)$ 
and strictly decreasing when evaluated at $u(t, \cdot)$. 
\section{Numerical experiments and viscosity} \label{s:ne}
We now go back to the numerical experiments in~\cite{ACT}, which as we have seen provide the wrong intuition concerning Question~\ref{q:q}. A possible explanation of the reason why the numerical evidence is not reliable is given by the following argument. 

The numerical results in~\cite{ACT} have been obtained by relying on a Lax-Friedrichs type scheme. The Lax-Friedrichs scheme is a finite volume scheme which is  very commonly used to construct numerical solutions of conservation laws, see~\cite{LeVeque} for an exhaustive discussion.  The Lax-Friedrichs scheme contains a large amount of what is called \emph{numerical viscosity}: very loosely speaking, the {numerical viscosity} is a collection of finite difference terms which are the numerical counterpart of some analytical viscosity, i.e. of some second order term. In other words, the presence of the numerical viscosity implies that the model equation for the Lax-Friedrichs scheme for the conservation law~\eqref{e:l} is actually the \emph{viscous} conservation law
\begin{equation} \label{e:vl}
          \partial_t u + \partial_x \big[ u V({ u}) \big] = \nu \partial^2_{xx} u,
\end{equation}
where the viscosity coefficent $\nu >0$ is of the same order of the space mesh, see~\cite{LeVeque}. 
When  the Lax-Friedrichs scheme is applied to the nonlocal conservation law~\eqref{e:nl}, the presence of the numerical viscosity implies that the model equation is
\begin{equation}
\label{e:vnl}
         \partial_t u + \partial_x \big[ u V({ u \ast \eta}) \big] = \nu \partial^2_{xx} u. 
\end{equation}
This in turn implies that in order to get some insight on the discrepancy between the numerical evidence in~\cite{ACT} and the analytic counterexamples in~\cite{CCS} it might be useful to consider the family of Cauchy problems\footnote{Existence and uniqueness results for the Cauchy problem~\eqref{e:vnlee} can be obtained 
by combining a fixed point argument with fairly standard parabolic estimates, see~\cite[\S2.1]{CCS}} 
\begin{equation}
\label{e:vnlee}
     \left\{ 
     \begin{array}{ll}
     \partial_t u_\ep + \partial_x \big[ u_\ep V({ u_\ep \ast \eta_\ep}) \big] = \nu \, \partial^2_{xx} u_\ep \\
     u_\ep (0, \cdot) = \bar u \\
    \end{array}
    \right.
\end{equation}
and pose the following ``viscous counterpart" of Question~\ref{q:q}. 
\begin{q}
\label{q:q2}
Does $u_\ep$, solution of~\eqref{e:vnlee}, converge to the solution of~\eqref{e:id},\eqref{e:vl} as $\ep~\to~0^+$? 
\end{q}
The answer to Question~\ref{q:q2} is largely positive, and it is given by the following result. 
\begin{theorem} \label{t:main}
Assume~\eqref{e:eta}, fix $\nu>0$ and $T>0$ and assume that the function $V: \R \to \R$ is Lipschitz continuous. If $\bar u \in L^1 (\R) \cap L^\infty (\R)$, then 
the solution of~\eqref{e:vnlee} converge to the solution of~\eqref{e:id},\eqref{e:vl} strongly in $L^2 ([0, T] \times \R)$ as 
$\ep~\to~0^+$. 
\end{theorem}
The proof of Theorem~\ref{t:main} is provided in~\cite{CCS} and applies in greater generality to the case of several space dimensions: we refer 
to~\cite[Theorem 1.1]{CCS} for the precise statement. Note furthermore that Theorem~\ref{t:main} was established in~\cite{CalderoniPulvirenti} under the additional assumptions that the initial datum $\bar u$ is regular and that $V(u)=u$. 

We can now consider the family of Cauchy problems~\eqref{e:vnlee}, keep the nonlocal parameter $\ep>0$ fixed, vary the viscosity parameter $\nu$ and 
consider the inviscid limit $\nu \to 0^+$. In this way we recover the inviscid nonlocal problem~\eqref{e:nlee}: more precisely, \cite[Proposition 1.2]{CCS} states that the solutions of~\eqref{e:vnlee} converge to the solution of~\eqref{e:nlee} when $\nu \to 0^+$. Finally, we recall that a celebrated result by 
Kru{\v{z}}kov~\cite{Kruzkov} states that the solutions of~\eqref{e:id},\eqref{e:vl} converge to the entropy admissible solution of~\eqref{e:id},\eqref{e:l} when $\nu \to 0^+$.

We now put together all the previous convergence results and we combine them with the counterexamples mentioned in~\S\ref{s:nltl}. We denote by $u_{\ep \nu}$ the solution of the viscous nonlocal equation at the first line of~\eqref{e:vnlee} to stress that it depends on both the nonlocal parameter $\ep$ and the viscosity parameter $\nu$. We arrive at the following diagram: 
\begin{equation*}
\minCDarrowwidth70pt
\begin{CD}
 \partial_t u_{\ep \nu} \! + \!\partial_x  \big[ u_{\ep \nu} V(u_{\ep \nu} \ast \eta_\ep) \big]\! =\! \nu \ \partial^2_{xx} u_{\ep \nu}   \!\!  @> \ep \to 0^+  >  \text{Theorem~\ref{t:main}}>  \! \! \partial_t u_\nu \! + \! \partial_x \big[ u_\nu V (u_\nu) \big] \! = \! \nu \ \partial^2_{xx}  u_\nu \\
@V \nu \to 0^+  V \text{ \cite[Proposition 1.2]{CCS}} 
V        @V \nu \to 0^+ V \text{Kru{\v{z}}kov~\cite{Kruzkov} } V\\
 \partial_t u_{\ep } + \partial_x \big[ u_{\ep } V(u_{\ep } \ast \eta_\ep) \big] = 0      @>  \ep \to 0^+ > \text{False}>   \partial_t u + \partial_x \big[ u V(u) \big] = 0
\end{CD}
\end{equation*}
We can now go back to the numerical evidence erroneously suggesting a positive answer to Question~\ref{q:q}. A possible explanation is the following: the numerical experiments were designed to test the convergence of~$u_\ep$ to the entropy admissible solution $u$. However, owing to the numerical viscosity, what the numerical experiments were actually testing was the convergence of~$u_{\ep \nu}$ to~$u_\nu$, which holds true owing to  Theorem~\ref{t:main}. In other words, the numerical schemes were designed to provide an answer to Question~\ref{q:q}, but as a matter of fact they provide an answer to Question~\ref{q:q2}. Since the two questions have opposite answers, the numerical schemes provide the wrong intuition.  
\begin{figure}[h!]
\begin{center}
\bigskip
\includegraphics[scale=0.2]{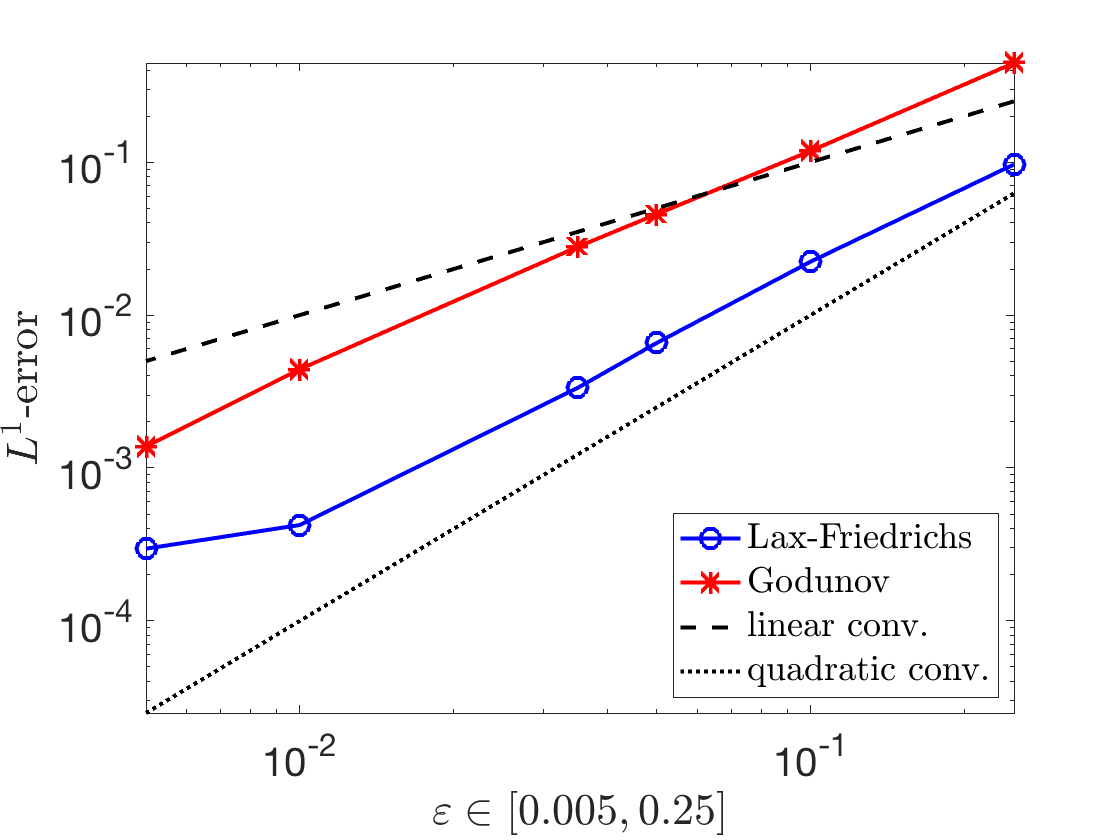}
\caption{$L^1$-error at $t=2$, for different values of $\ep$, comparing the solution of~\eqref{e:nlee} to the entropy admissible solution of~\eqref{e:id},\eqref{e:l} computed with Lax-Friedrichs and Godunov type schemes in the case where~\eqref{e:nlee} is the same as in~\cite[\S5.1]{CCS}. The space mesh is fixed and it is $h=0.001$} \label{f:1}
\end{center}
\end{figure}
This explanation is validated by recent numerical experiments collected in~\cite{CCGS}. In particular, in~\cite{CCGS}, we have used the Lax-Frierichs type scheme to test the nonlocal-to-local limit from~\eqref{e:nlee} to~\eqref{e:id},\eqref{e:l} in the case of the counterexamples mentioned in~\S\ref{s:nltl}. 
More precisely, we have computed the numerical solution of~\eqref{e:nlee} in the case where~\eqref{e:nlee} is the same as in the counterexamples. Next, we have computed the $L^1$ norm  (evaluated at a given positive time $t>0$)  of the difference between the numerical solution of~\eqref{e:nlee} and the numerical entropy admissible solution of~\eqref{e:id},\eqref{e:l}. In Figure~\ref{f:1} we display some of the results concerning one of the counterexamples, more precisely the one discussed in~\cite[\S5.1]{CCS}. 
\begin{figure}[h!]
\begin{center}
\bigskip
\includegraphics[scale=0.15]{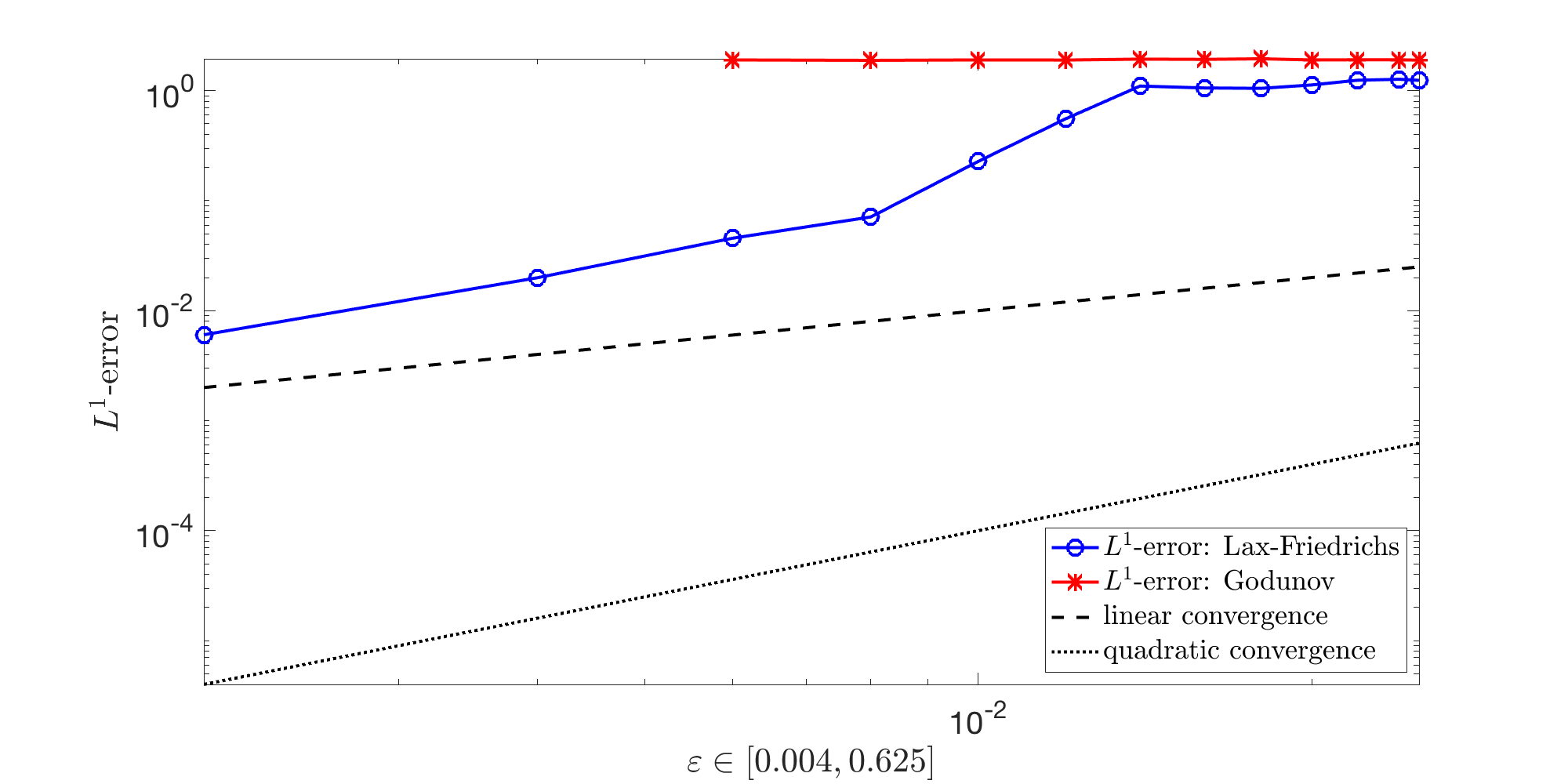}
\caption{$L^1$-error at $t=2$, for different values of $\ep$, between the solutions of the nonlocal equations~\eqref{e:nlee} and the entropy solution of~\eqref{e:id},\eqref{e:l} computed with Godunov and Lax Friedrichs schemes in the case where~\eqref{e:nlee} is the same as in~\cite[\S5.2]{CCS}. The space mesh $h$ depends on $\ep$ and the relation is $\ep = 1000 h^2$. The $L^1$ error of the Godunov scheme is much larger for small values of $\ep$.} \label{f:2}
\end{center}
\end{figure}
The blue line refers to the $L^1$ error between the numerical solutions obtained by the Lax-Frierichs type scheme and strongly suggests that the $L^1$ error is converging to $0$ as $\ep \to 0^+$, i.e. it erroneously suggests a positive answer to Question~\ref{q:q}. The red line refers to the $L^1$ error between the numerical solutions obtained by a Godunov type scheme. Godunov type schemes for the nonlocal conservation law~\eqref{e:nl} were introduced in~\cite{ChiarelloGoatin,FKG} and the reason why we used them to test the nonlocal-to-local limit is because the Godunov scheme is known to have a smaller amount of numerical viscosity than the Lax-Friedrichs scheme, see~\cite{Tadmor}.  In the example studied in Figure~\ref{f:1}, the numerical results obtained with both the Godunov and the Lax-Friedrichs scheme erroneously suggest convergence in the nonlocal-to-local limit. However, in other cases there is a difference between the two schemes. For instance, Figure~\ref{f:2} displays some of the numerical results concerning the counterexample discussed in~\cite[\S5.2]{CCS}: there is a remarkable difference between the Lax-Friedrichs and the Godunov scheme. Indeed, the numerical results obtained with the Lax-Friedrichs type scheme erroneously suggest convergence, whereas the numerical results  obtained with the Godunov type scheme are more consistent with the analytic results, which rule out convergence. This is consistent with the fact that the Godunov scheme contains less numerical viscosity than the Lax-Friedrichs scheme, see~\cite{Tadmor}. 
\section{Anisotropic traffic models: total variation blow up and open questions} \label{s:ak}
In recent years, several authors have been focusing on~\eqref{e:nl} in the case where the function $V$ is decreasing, $V'<0$, the initial datum $\bar u$ is nonnegative and 
the convolution kernel $\eta$ in equation~\eqref{e:l} is completely anisotropic (i.e., it is supported on $]-\infty, 0]$). This case is extremely relevant for the applications to vehicular traffic models. 
Indeed, in these models $u$ represents the density of cars (and is therefore nonnegative) and $V$ their speed. The function $V$ is evaluated at $u \ast \eta$ because the model postulates that drivers regulate their speed based on the density of cars around them. 
 The fact that the function $V$ is decreasing is a classical assumption in traffic models and takes into account the fact that drivers tend to slow down when the traffic is congested, and conversely to speed up when the traffic is light. If the convolution kernel is supported on the interval $]-\infty, 0]$, then the convolution kernel $u\ast \eta$ evaluated at the point $x$ only depends on the value of $u$ on the interval $[x, + \infty[$. 
In other words, choosing an anisotropic convolution kernel aims at modeling  the fact that drivers only look forward, not backward, and hence their speed only depends on the downstream traffic density.   
 
To avoid some technicalities, in the following  we focus on the case 
\begin{equation}
\label{e:traffic}
     V(u) = 1-u, \quad \eta= \mathbbm{1}_{[-1, 0]}, \quad 0 \leq \bar u \leq 1,
\end{equation} 
but as a matter of fact the following  discussion applies to more general cases than~\eqref{e:traffic}.
In the previous formula, $\mathbbm{1}_{[-1, 0]}$ denotes the characteristic function of the interval $[-1, 0]$.  Note that, strictly speaking, the regularity assumptions on the function $\eta$ given in~\eqref{e:eta} are violated when $\eta= \mathbbm{1}_{[-1, 0]}$. Notwithstanding the lack of regularity, in~\cite{BlandinGoatin} the authors established existence and uniqueness results for the Cauchy problem~\eqref{e:nl},\eqref{e:id} and, remarkably, by exploiting the anisotropy of the kernel established better a-priori estimates on the solution than those available in the smooth case~\eqref{e:eta}.  In particular, they established a maximum principle: under~\eqref{e:traffic}, the solution of~\eqref{e:nl},\eqref{e:id} satisfies $0 \leq  u  \leq 1$. 
To complete the picture, we point out that the counterexamples exhibited in~\cite{CCS} \emph{do not} apply in the case~\eqref{e:traffic}. 

Summing up, the case~\eqref{e:traffic} is very relevant from the modeling viewpoint, stronger analytic results apply and the counterexamples do not work. This yields the following refinement of Question~\ref{q:q}. 
\begin{q}\label{q:q3} 
Does $u_\ep$, solution of~\eqref{e:nlee}, converge to the entropy admissible solution of~\eqref{e:id},\eqref{e:l} as $\ep \to 0^+$, provided~\eqref{e:traffic} holds true? 
\end{q}
Question~\ref{q:q} is presently open and it is the object of current investigation. However, some progress have been recently achieved in~\cite{CCS2}. Before discussing the results in~\cite{CCS2}, we need some preliminary considerations. 

Assume~\eqref{e:traffic}, then, owing to the maximum principle, the solution of the Cauchy problem~\eqref{e:nlee} satisfies
the uniform bound 
$$ 
    \| u_\ep \|_{L^\infty} \leq 1, \quad \text{for every $\ep>0$}.
 $$   
 This yields compactness in the weak-$^\ast$ topology and implies that we can extract a subsequence that converges to some limit function $w$ weakly-$^\ast$ in $L^\infty (\R^+ \times \R)$. Note however that, owing to the nonlinear nature of the problem, 
nothing a priori tells us that the limit $w$ is a distributional solution (let alone entropy admissible) of the conservation law~\eqref{e:id},\eqref{e:l}. 
A natural strategy to establish a positive answer to Question~\ref{q:q2} is hence to look for compactness in some \emph{strong} topology. 
A fairly classical argument to establish strong $L^1$ compactness combines the Helly-Kolmogorov  Compactness Theorem with a uniform bound on the total variation, i.e. an estimate like 
\begin{equation}
\label{e:tv}
      \mathrm{TotVar} \, u_\ep (t, \cdot) \leq C, \quad \text{for every $t>0, \, \ep>0$ and for some constant $C>0$}.
\end{equation}
This yields the following question: 
\begin{q}\label{q:q4} 
Assume~\eqref{e:traffic} and that $\mathrm{TotVar} \,  \bar u$ is finite. Does $u_\ep$, solution of~\eqref{e:nlee}, satisfy the uniform bound~\eqref{e:tv}?
\end{q}
Before addressing Question~\ref{q:q4} we make some preliminary remarks. 
First,  the semigroup of entropy admissible solutions of~\eqref{e:id},\eqref{e:l} is total variation decreasing, i.e. 
\begin{equation}
\label{e:tvd}
      \mathrm{TotVar} \, u (t, \cdot) \leq \mathrm{TotVar} \,  \bar u, \quad \text{for every $t>0$},
\end{equation}
provided  $\mathrm{TotVar} \,  \bar u$ is finite. In other words, the entropy admissible solution of~\eqref{e:id},\eqref{e:l} satisfies estimate~\eqref{e:tv} with $C= \mathrm{TotVar} \,  \bar u$.
Second, numerical experiments discussed in~\cite{BlandinGoatin} suggest that, under~\eqref{e:traffic}, the semigroup of solutions of~\eqref{e:nlee} is also total variation decreasing, and hence in particular that the answer to Question~\ref{q:q4} is positive. 
Third, by combining the maximum principle with the monotonicity preserving property established in~\cite{BlandinGoatin} one can show that that, under~\eqref{e:traffic}, if the initial datum $\bar u$ is \emph{monotone}, then the total variation does not increase in time, i.e.~\eqref{e:tv} is satisfied with $C= \mathrm{TotVar} \,  \bar u$. In other words, we know that the answer to Question~\ref{q:q4} is positive provided the initial datum is monotone.

Notwithstanding  the numerical evidence and the positive answer in the case of monotone data, a counterexample constructed in~\cite{CCS2} shows that the answer to Question~\ref{q:q4} is in general negative. More precisely, there is an initial datum $\bar u$ such that $\mathrm{TotVar} \,  \bar u$ is finite and the solution of the Cauchy problem~\eqref{e:nlee} 
 satisfies 
$$
    \sup_{\ep >0}  \mathrm{TotVar} \, u_\ep (t, \cdot) = + \infty, \; \text{for every $t>0$},
$$
which in particular implies that~\eqref{e:tv} cannot be true. 

The fact that the answer to Question~\ref{q:q4} is negative does not by any mean imply that the answer to Question~\ref{q:q3} is also negative. 
However, it rules out the most classical and natural strategy to achieve an hypothetical positive answer to Question~\ref{q:q3}. 
Note, furthermore, that the initial datum $\bar u$ in~\cite{CCS2}  has finite total variation and attains values in the physical range $0 \leq \bar u \leq 1$, but it is also highly oscillating and it is unlikely to describe a realistic initial density of vehicles in some real-word applications. In principle it might be possible that, under~\eqref{e:traffic}, the uniform bound~\eqref{e:tv} holds true provided $\bar u$ is an initial datum with finite total variation which satisfies some further condition making it more ``realistic". Even if this were true, however, the counterexample in~\cite{CCS2} would provide some useful information because it implies that~\eqref{e:tv} cannot be established by relying only on the maximum principle and on the boundedness of $\mathrm{TotVar} \,  \bar u$. To establish~\eqref{e:tv} in the case of ``realistic" initial data one should likely rely on some more refined information on the structure of the solution, which  is in general harder to obtain.

\section*{Acknowledgments} GC is partially supported by the Swiss National Science 
Foundation grant 200020-156112 and by the
ERC Starting Grant 676675 FLIRT. MG was partially
supported by the Swiss National Science Foundation grant P300P2-167681. LVS is a member of the GNAMPA
group of INDAM and of the PRIN National Project ``Hyperbolic Systems of Conservation
Laws and Fluid Dynamics: Analysis and Applications''. 
\bibliographystyle{plain}
\bibliography{hyp2018}
\end{document}